\documentclass[9pt]{article}
\input epsf.tex
\usepackage{latexsym}
\usepackage{amsfonts,amssymb}
\usepackage{epsfig}

\begin{document}
\title{Compressing totally geodesic surfaces}
\author{Christopher J. Leininger}
\maketitle
\begin{abstract}
In this paper we prove that one can find surgeries arbitrarily close to infinity in the Dehn surgery space of the figure eight knot complement for which some immersed totally geodesic surface compresses.
\end{abstract}


\section{Introduction}

Suppose that $M$ is a compact, orientable, irreducible $3$-manifold with torus boundary $\partial M \cong T^{2}$.
Recall that a {\em slope} on $\partial M$ is an isotopy class of simple closed curves on $\partial M$.
The {\em distance} between a pair of slopes $\alpha$ and $\beta$, denoted $\Delta (\alpha , \beta)$, is the absolute value of their algebraic intersection number.

Given a slope $\alpha$, let $M(\alpha)$ denote the manifold obtained from $\alpha$-Dehn surgery on $\partial M$.
It was shown in \cite{CGLS} that if $F$ is a closed, orientable, embedded, incompressible surface in $M$ admitting no incompressible annulus with one boundary component contained in $F$ and the other in $\partial M$, and $F$ compresses in $M(\alpha)$ and $M(\beta)$, then $\Delta (\alpha,\beta) \leq 2$.
This result was later improved to $\Delta (\alpha,\beta) \leq 1$ \cite{Wu}.

In this paper we study the analog of the above result for immersed surfaces.
When the interior of $M$ admits a hyperbolic structure of finite volume, closed, immersed, totally geodesic surfaces are incompressible surfaces having no incompressible annulus as above.
These surfaces have good incompressibility properties, for they remain incompressible in all but a finite number of surgeries on one cusp (this can be deduced from \cite{Ma} and \cite{T}, see also \cite{Bt} for explicit bounds).

A well known example of a complete hyperbolic manifold with one cusp is the complement of the figure eight knot in $S^{3}$ \cite{Ri}, which we denote by $M_{8}$.
Using arithmeticity of $M_{8}$, one can show that it contains infinitely many closed, immersed totally geodesic surfaces \cite{MR2},\cite{R2}.

Our main result is\\

\noindent
{\bf Theorem 5.1.1 } {\em Suppose $4 | p$ and $ 3 \nmid p$.
Then, for any $q$ with $gcd(p,q)=1$ there exists infinitely many noncommensurable, closed, orientable, immersed, totally geodesic surfaces in $M_{8}$ which compress in $M_{8}(\frac{p}{q})$.}\\

Theorem 5.1.1 implies that there is no global bound on the number of surgeries which one must omit to guarantee incompressiblity of all closed totally geodesic surfaces.
More precisely, we have\\

\noindent
{\bf Corollary } {\em There are infinitely many surgeries on the cusp of $M_{8}$ such that for each surgery, some immersed, closed, totally geodesic surface compresses.}\\

As $M_{8}$ is the interior of a compact $3$-manifold with torus boundary, and since every surface in $M_{8}$ compresses in $M_{8}(\frac{1}{0})$, Theorem 5.1.1 implies the following result, indicating there is no analog in the immersed settings of \cite{CGLS} and \cite{Wu}.\\

\noindent
{\bf Theorem 4.4.1 } {\em There exists a compact, orientable, irreducible $3$-manifold $M$, with torus boundary, having the following property.
Given any positive integer, $n$, there exists a closed, orientable, immersed, incompressible surface $F \looparrowright M$ with no incompressible annulus joining $F$ and $\partial M$, such that $F$ compresses in $M(\alpha)$ and $M(\beta)$ and $\Delta(\alpha,\beta) > n$.}\\

The paper is organized as follows:
We briefly discuss general $3$-manifold topology in Section 2.
In Section 3, we describe a few of the necessary results from hyperbolic geometry, and Section 4 contains results from the theory of arithmetic hyperbolic manifolds needed to prove Theorem 5.1.1.
We complete this section by proving Theorem 4.4.1, assuming Theorem 5.1.1.
In Section 5, we prove Theorem 5.1.1 and describe some other examples obtainable using a similar construction.\\

{\bf Acknowledgments}:  I would like to thank my advisor Alan Reid for many useful conversations and suggestions.

\section{$3$-manifolds and incompressible surfaces}

In this section we collect some of the basic facts and definitions from $3$-manifold topology, see \cite{H} and \cite{Ro} for more details.

Let $M$ denote a compact, orientable, irreducible $3$-manifold with torus boundary $\partial M \cong T^{2}$.
Slopes on $\partial M$ are in a $2$ to $1$ correspondence with primitive elements of $\pi_{1}(\partial M) \cong H_{1}(\partial M,{\mathbb Z}) \cong {\mathbb Z} \oplus {\mathbb Z}$, with the ambiguity coming from the lack of orientations.
We will often ignore this ambiguity, making no distinction between slopes and primitive elements of $\pi_{1}(\partial M)$.

Let $\mu$ and $\lambda$ denote generators for $\pi_{1}(\partial M)$.
Slopes on $\partial M$ then correspond to elements $r \in {\mathbb Q} \cup \{ \infty \}$ by
associating $\mu^{p} \lambda^{q}$ with $r = \frac{p}{q}$, where $\frac{p}{q}$ is in lowest terms and $\infty = \frac{1}{0}$.
Accordingly, we will denote the slope $\alpha$ by its associated $r \in {\mathbb Q} \cup \{ \infty \}$. 
One checks that if $\alpha = \frac{p}{q}$ and $\beta = \frac{p'}{q'}$ then $\Delta(\alpha,\beta) = |pq' - p'q|$.

Given a slope $\alpha$ on $\partial M$, one can form a new manifold $M(\alpha)$ by {\em $\alpha$-Dehn surgery} as follows.
Let $S^{1} \times D^{2}$ be a solid torus.
Choosing a homeomorphism $h : \partial (S^{1} \times D^{2}) \rightarrow \partial M$, so that $h(* \times \partial D^{2})$ represents $\alpha$, we can glue $S^{1} \times D^{2}$ to $M$ by identifying points $x$ and $h(x)$.
The resulting space is a manifold, and up to homeomorphism, depends only on $\alpha$.
Note that there is a natural inclusion $i : M \rightarrow M(\alpha)$. 

Let $F \ncong S^{2},D^{2}$ be an orientable surface and $f : F \rightarrow N$ an immersion into an orientable $3$-manifold $N$.
We will say that $(f,F)$ is an {\em incompressible surface} if $f_{*}:\pi_{1}(F) \rightarrow \pi_{1}(N)$ is injective, and {\em compressible} otherwise.
We may often refer to an incompressible or compressible surface $F$, with the mapping implicit.
By Dehn's Lemma and the Loop Theorem, if $f$ is a proper embedding, these definitions agrees with the usual ones \cite{H}.

It will often be the case that the $3$-manifold $M$ we are considering is the interior of a compact $3$-manifold $\overline{M}$ with torus boundary.
In this case, $M(\alpha)$ will be the manifold $\overline{M}(\alpha)$, where $\alpha$ denotes both the slope on $\partial \overline{M}$ and the associated homotopy class of curves in the end of $M$.
A homotopy class of curves in $M$ which have a representative in $\partial \overline{M}$ is called {\em peripheral}.
We will say that a loop in $M$ is peripheral if its homotopy class is.

\section{Hyperbolic geometry}

In this section we review a few of the basics of hyperbolic geometry, for more details, see \cite{BP},\cite{Ra}, or \cite{T}.

\subsection{Hyperbolic space and its isometries}

We will use the upper half space model for hyperbolic $3$-space, ${\mathbb H}^{3}$.
That is,
$${\mathbb H}^{3} = \{ (z,t) \, | \, z = x + i y \in {\mathbb C} \, , \, t > 0 \}$$
with the complete Riemannian metric
$$ds^{2} = \frac{dx^{2} + dy^{2} + dt^{2}}{t^{2}}.$$
${\mathbb H}^{3}$ is compactified by $S^{2}_{\infty} = \widehat{\mathbb C}$ and all orientation preserving isometries are conformal extensions of conformal maps of $\widehat{\mathbb C}$.
Thus, $PSL_{2}({\mathbb C})$ is the full group of orientation preserving isometries of ${\mathbb H}^{3}$, acting by extension of linear fractional transformations on $\widehat{\mathbb C}$.

Similarly, the upper half plane model of the hyperbolic plane is
$${\mathbb H}^{2} = \{ (x,t) \, | \, x \in {\mathbb R} \, , \, t > 0 \}$$
with metric
$$ds^{2} = \frac{dx^{2} + dt^{2}}{t^{2}}.$$
${\mathbb H}^{2}$ is naturally compactified by $S^{1}_{\infty}=\widehat{\mathbb R}$.
$PSL_{2}({\mathbb R})$ is the full group of orientation preserving isometries of ${\mathbb H}^{2}$ acting by extension of linear fractional transformations on $\widehat{\mathbb R}$.
Note that we can view ${\mathbb H}^{2}$ as a submanifold of ${\mathbb H}^{3}$ embedded totally geodesically.
When this is done, the action of $PSL_{2}({\mathbb R})$ on ${\mathbb H}^{2}$ is the restriction to ${\mathbb H}^{2}$ of the action of $PSL_{2}({\mathbb R}) \subset PSL_{2}({\mathbb C})$ on ${\mathbb H}^{3}$.
When we identify ${\mathbb H}^{2}$ with the upper half plane of ${\mathbb C}$ (in the obvious way), the action of $PSL_{2}({\mathbb R})$ is by linear fractional transformations.\\

\noindent
{\bf Remark: } Let $P : SL_{2}({\mathbb C}) \rightarrow PSL_{2}({\mathbb C})$ denote the quotient map.
Whenever we refer to a matrix for an element $g \in PSL_{2}({\mathbb C})$, we mean a matrix $A \in P^{-1}(g)$.\\

Given $\gamma \in PSL_{2}({\mathbb C})$, the trace of $\gamma$, $Tr(\gamma)$, is well defined up to sign.
We say that $\gamma$ is {\em elliptic} if $Tr(\gamma) \in (-2,2)$, {\em parabolic} if $Tr(\gamma) = \pm 2$, and {\em hyperbolic} otherwise.
If $\gamma \neq 1$, then $Tr(\gamma)$ is a complete invariant of the conjugacy class of $\gamma$.
See \cite{T} for a geometric description of the action on ${\mathbb H}^{n}$, $n=2,3$.

\subsection{Kleinian groups}

A discrete subgroup $\Gamma \subset PSL_{2}({\mathbb C})$ is called a {\em Kleinian group}.
Discreteness of $\Gamma$ is equivalent to the action on ${\mathbb H}^{3}$ being properly discontinuous.
Proper discontinuity easily implies $\Gamma$ is torsion free if and only if it contains no elliptic elements.
$\Gamma$ is said to be {\em elementary} if it contains an abelian subgroup of finite index, and {\em non-elementary} otherwise.

For the remainder of this section, let $\Gamma$ represent a torsion free Kleinian group.
We let $M_{\Gamma} = {\mathbb H}^{3} / \Gamma $ denote the quotient {\em hyperbolic 3-manifold} (with its induced metric), and let
$$p : {\mathbb H}^{3} \rightarrow M_{\Gamma}$$
denote the projection.
Note that $({\mathbb H}^{3},p)$ is the universal covering of $M_{\Gamma}$, and $\Gamma$ is the group of covering transformations, so that $\pi_{1}(M_{\Gamma}) \cong \Gamma$.
When we wish to refer to this isomorphism explicitly, we will write it as
$$ \psi : \pi_{1}(M_{\Gamma}) \rightarrow \Gamma.$$

We say that $\Gamma$ has {\em finite co-volume} (resp. is {\em co-compact}) if $M_{\Gamma}$ has finite total volume (resp. is compact).
A {\em cusp} of $M_{\Gamma}$ is a subset of $M_{\Gamma}$ isometric to a set of the form $B / \Gamma_{P}$ where
$$B = \{ (x,y,t) \in {\mathbb H}^{3} \, | \, t > 1 \}$$
and $\Gamma_{P} \cong {\mathbb Z} \oplus {\mathbb Z}$ is a subgroup of $PSL_{2}({\mathbb C})$ consisting entirely of parabolics which stabilizes $B$.
It can be shown (see \cite{BP} for example) that when $\Gamma$ has finite co-volume, $M_{\Gamma}$ is the interior of a compact manifold with toroidal boundary, and that every boundary component of $M_{\Gamma}$ has a product neighborhood whose intersection with $M_{\Gamma}$ consists of cusps of $M_{\Gamma}$.

Let $\gamma$ be a free homotopy class of essential loops in $M_{\Gamma}$.
We say that $\gamma$ is hyperbolic (resp. parabolic) if $\psi(\gamma')$ is hyperbolic (resp. parabolic) where $\gamma' \in \pi_{1}(M_{\Gamma})$ is a representative of the conjugacy class determined by $\gamma$. 
If $\gamma$ is hyperbolic, then there exists a unique, length minimizing, geodesic representative for $\gamma$ given by $Ax(\psi(\gamma')) / <\psi(\gamma')>$, where $Ax(\psi(\gamma'))$ is the axis for $\psi(\gamma')$.
If $\gamma$ is parabolic, then there exists a sequence of loops $\gamma_{n}$ representing $\gamma$ such that the length of $\gamma_{n}$ goes to $0$ as $n$ goes to infinity.

As a homotopy class of loops is peripheral if and only if it has representatives lying entirely in a cusp of $M_{\Gamma}$, the following is straightforward.\\

\noindent
{\bf Theorem 3.2.1 } {\em Let $M_{\Gamma}$ be a finite volume hyperbolic $3$-manifold.
A homotopy class of essential loops in $M_{\Gamma}$ is peripheral if and only if it is parabolic.}\\

\noindent
{\bf Notation: } It is common to blur the distinction between a particular cusp and the end of the manifold corresponding to that cusp.
We will follow this convention, referring to both objects as cusps.
The context will make it clear which we are referring to.

\subsection{Fuchsian groups}

Given any subgroup $G \subset PSL_{2}({\mathbb C})$ and any (geometric) circle ${\cal C} \subset \widehat{{\mathbb C}}$ (i.e. circle or line in ${\mathbb C}$), define
$$Stab_{G}({\cal C}) = \{ g \in G \, | \, g({\cal C}) = {\cal C} , \, g \mbox{ preserves the components of } \widehat{\mathbb C} \setminus {\cal C} \}$$
A {\em Fuchsian group} is defined to be a discrete subgroup of $Stab_{PSL_{2}({\mathbb C})}({\cal C})$, for some circle ${\cal C}$.
For any circle ${\cal C} \subset \widehat{\mathbb C}$ there exists $g \in PSL_{2}({\mathbb C})$ such that $g({\cal C}) = \widehat{\mathbb R}$.
It follows that $g Stab_{PSL_{2}({\mathbb C})}({\cal C}) g^{-1} = PSL_{2}({\mathbb R}) $.
Thus, a Fuchsian group is a Kleinian group conjugate into $PSL_{2}({\mathbb R})$ by an element $g \in PSL_{2}({\mathbb C})$.

Any circle ${\cal C}$ on $\widehat{{\mathbb C}}$ is the boundary of a hyperbolic plane $P_{\cal C} \cong {\mathbb H}^{2}$ embedded totally geodesically in ${\mathbb H}^{3}$ and conversely.
In the notation of the previous paragraph, we have $P_{\cal C} = g^{-1}({\mathbb H}^{2})$, where we view ${\mathbb H}^{2} \subset {\mathbb H}^{3}$ as in Section 3.1.
If $\Gamma$ is a torsion free Fuchsian group stabilizing ${\cal C}$ (hence $P_{{\cal C}}$), $S_{\Gamma} = P_{\cal C} / \Gamma$ is a {\em hyperbolic surface}, and $\pi_{1}(S_{\Gamma}) \cong \Gamma$.
$\Gamma$ has {\em finite co-area} (resp. is {\em co-compact}) if $S_{\Gamma}$ has finite total area (resp. is compact).

Suppose now that $\Gamma$ is a finite co-volume torsion free Kleinian group such that there exists a circle ${\cal C} \subset \widehat{\mathbb C}$ for which $\Gamma' = Stab_{\Gamma}({\cal C})$ has finite co-area.
This induces a proper, totally geodesic immersion
$$S_{\Gamma'} \looparrowright M_{\Gamma}$$
It is immediate that any such surface is incompressible.
Since the only (complete) totally geodesic surfaces in ${\mathbb H}^{3}$ are hyperbolic planes, any proper, totally geodesic immersion of an orientable, hyperbolic surface into $M_{\Gamma}$ factors through an immersion of this type.
That is, if $f : F \looparrowright M_{\Gamma}$ is a proper, totally geodesic immersion of an orientable, hyperbolic surface, then with the notation above, we have that $\psi \circ f_{*}(\pi_{1}(F)) \subset Stab_{\Gamma}({\cal C})$ for some circle ${\cal C} \subset \widehat{\mathbb C}$.

We will be primarily interested in closed, orientable, immersed, totally geodesic surfaces.
The classification of isometries of ${\mathbb H}^{n}$, $n = 2,3$, easily implies\\

\noindent
{\bf Theorem 3.3.1 } {\em If $F$ is a closed, orientable, immersed, totally geodesic surface in a finite volume hyperbolic $3$-manifold $M_{\Gamma}$, then every free homotopy class of essential loops in $F$ is hyperbolic.}\\

Theorem 3.2.1 and Theorem 3.3.1 together imply\\

\noindent
{\bf Corollary 3.3.2 } {\em If $F$ is a closed, orientable, immersed, totally geodesic surface in a finite volume hyperbolic $3$-manifold $M_{\Gamma}$, then there are no free homotopy classes of essential loops in $F$ that are peripheral.}\\

\noindent
{\bf Remark: } The surfaces $S_{\Gamma'}$ corresponding to Fuchsian subgroups $Stab_{\Gamma}({\cal C}) = \Gamma' \subset \Gamma$ are orientable, although it may be that the map $S_{\Gamma'} \looparrowright M_{\Gamma}$ factors through a covering of a non-orientable totally geodesic surface $\Sigma \looparrowright M_{\Gamma}$.
This will happen if and only if there is a $g \in \Gamma$ such that $g({\cal C}) = {\cal C}$ and $g$ {\em exchanges} the components of $\widehat{\mathbb C} \setminus {\cal C}$.

\section{Arithmetic manifolds}

In this section, we discuss background and results from the theory of arithmetic manifolds, applications to the figure eight knot, and we give a proof of Theorem 4.4.1.
All theorems in this section (with the exception of Theorem 4.4.1) are known and a few proofs have been included for the sake of completeness.
For more details, see \cite{MR1}.

\subsection{Arithmetic Fuchsian Groups}

Let $A$ be a quaternion algebra over ${\mathbb Q}$ with Hilbert symbol $\left( \frac{a,b}{\mathbb Q} \right)$.
That is, $A$ is a $4$-dimensional algebra over ${\mathbb Q}$ having basis $1,{\bf i},{\bf j},{\bf k}$ with multiplication defined so that $1$ is a multiplicative identity, and 
$$ {\bf i}^{2} = a \cdot 1 \, , \, {\bf j}^{2} = b \cdot 1 \, , \, {\bf ij}=-{\bf ji}={\bf k}$$
where $a,b \in {\mathbb Q}^{*}$.
$A$ admits an anti-involution $x \longmapsto \overline{x}$ called {\em conjugation}.
That is, $x \longmapsto \overline{x}$ is an involution of the vector space and $\overline{x \cdot y} = \overline{y} \cdot \overline{x}$.
Conjugation is given by
$$\overline{x} = \overline{x_{0} + x_{1} {\bf i} + x_{2} {\bf j} + x_{3} {\bf k}} = x_{0} - x_{1} {\bf i} - x_{2} {\bf j} - x_{3} {\bf k} $$

The {\em (reduced) norm} and {\em (reduced) trace} of $x \in A$ are defined by $n(x) = x \overline{x}$ and $tr(x) = x + \overline{x}$ respectively.
We will view $tr$ and $n$ both as maps to ${\mathbb Q}$.

We note that the quaternion algebra $A = \left( \frac{a,b}{{\mathbb Q}} \right)$ embeds into $M_{2}({\mathbb Q}(\sqrt{a}))$ by
$$\rho(x) = \rho(x_{0} + x_{1}{\bf i} + x_{2}{\bf j} + x_{3}{\bf k})=
\left( \begin{array}{cc}
x_{0}+x_{1} \sqrt{a} & b(x_{2}+x_{3} \sqrt{a})\\
x_{2}-x_{3} \sqrt{a} & x_{0} - x_{1} \sqrt{a}\\ \end{array} \right)$$
With this embedding, we have that $tr(x) = Tr(\rho(x))$ and $n(x) = det (\rho(x))$ where $Tr$ and $det$ are the usual trace and determinant of a square matrix, respectively.

An {\em order} ${\cal O}$ in a quaternion algebra $A$ over ${\mathbb Q}$ is finitely generated ${\mathbb Z}$-module contained in $A$ such that ${\cal O}$ spans $A$ over ${\mathbb Q}$ (${\cal O} \otimes_{\mathbb Z} {\mathbb Q} = A$) and ${\cal O}$ is a ring with $1$.
Given an order ${\cal O}$, define $ {\cal O}^{1} =  \{ x \in {\cal O} \, | \, n(x) = 1 \}$.

A useful result is the following \cite{MR1}.\\

\noindent
{\bf Theorem 4.1.1 } {\em If $A$ is a quaternion algebra over ${\mathbb Q}$ and ${\cal O}$ is an order in $A$, then $P \rho {\cal O}^{1}$ is a finite co-area Fuchsian group.
Moreover, $P \rho {\cal O}^{1}$ is co-compact if and only if $A$ is a division algebra.}\\

From this theorem, we make the following definition.
A Fuchsian group ${\Gamma}$ is said to be {\em derived from a quaternion algebra} (defined over ${\mathbb Q}$) if $\Gamma$ is conjugate into a subgroup of $P\rho{\cal O}^{1}$ of finite index, for some $A$ and ${\cal O}$ as above.

Recall that if $G$ is any group and $K$ and $H$ are subgroups, then $K$ and $H$ are said to be {\em commensurable} in $G$ if and only if there exists an element $g \in G$ such that $g H g^{-1} \cap K$ is a finite index subgroup of both $g H g^{-1}$ and $K$.
Any Fuchsian group commensurable in $PSL_{2}({\mathbb C})$ with a group derived from a quaternion algebra is said to be {\em arithmetic}.

\subsection{Bianchi groups and Fuchsian subgroups}

Let $d$ be a positive square-free integer, and let $k_{d}={\mathbb Q}(\sqrt{-d})$.
Let ${\cal O}_{d}$ be the ring of integers in $k_{d}$.
That is,
$${\cal O}_{d} = \{ a + b \sqrt{-d} \, | \, a,b \in {\mathbb Z} \} \mbox{ if } d \equiv 1 \mbox{ or } 2 \mbox{ (mod 4) }$$
and
$${\cal O}_{d} = \{ \frac{a + b \sqrt{-d}}{2} \, | \, a,b \in {\mathbb Z} \, , \, a \equiv b \mbox{ (mod 2) } \} \mbox{ if } d \equiv 3 \mbox{ (mod 4) }$$
A {\em Bianchi group} is any group of the form $PSL_{2}({\cal O}_{d})$ for some square-free $d \in {\mathbb Z}^{+}$. 
It is well known that the Bianchi groups all have finite co-volume and are non-co-compact (see \cite{MR1}).

We will be interested in Fuchsian subgroups of the Bianchi groups.
We now describe a certain integral invariant of such subgroups that we will be useful for us.
Any circle ${\cal C}$ in $\widehat{\mathbb C}$ can be described by a triple $(a,B,c)$ as the set of $z \in {\mathbb C}$ (and possibly the point $\infty$) such that
$$a|z|^{2} + B z + \overline{Bz} + c = 0$$
If $Stab_{PSL_{2}({\cal O}_{d})}({\cal C})$ is non-elementary, it can be shown (see \cite{M}) that we may choose $a,c \in {\mathbb Z}$ and $B \in {\cal O}_{d}$.
In this case, we will write $B = \frac{1}{2}(b_{1}+b_{2} \sqrt{-d})$ where $b_{1},b_{2} \in {\mathbb Z}$ , $b_{1} \equiv b_{2}$ (mod 2), and both congruent to $0$ (mod 2) if $d \equiv 1,2$ (mod 4).
We say that the triple $(a,B,c)$ is {\em primitive} if 
$$ gcd(a,\frac{b_{1}}{2},\frac{b_{2}}{2},c)=1 \mbox{ for } b_{1} \equiv b_{2} \equiv 0 \mbox{ (mod 2) }$$
and
$$ gcd(a,b_{1},b_{2},c)=1 \mbox{ for } b_{1} \equiv b_{2} \not \equiv 0 \mbox{ (mod 2) }$$
A primitive triple for a circle is unique up to sign.
If $(a,B,c)$ is a primitive triple for ${\cal C}$, we define the {\em discriminant} of ${\cal C}$ to be ${\cal D}({\cal C}) = |B|^{2} - ac$ (note that ${\cal D}({\cal C}) > 0$).

We denote the set of circles represented by primitive triples in ${\cal O}_{d}$ by $\Sigma_{d}$.
$PSL_{2}({\cal O}_{d})$ acts on $\Sigma_{d}$ by $T \cdot {\cal C} = T({\cal C})$.
That is, the M$\ddot{o}$bius transformation represented by $T \in PSL_{2}({\cal O}_{d})$ takes the circle ${\cal C} \in \Sigma_{d}$ to the circle $T({\cal C})$, and $T({\cal C})$ is represented by a primitive triple.

Now define
$$ {\cal H}_{d} = \left\{ \left( \begin{array}{cc}
a & B \\
\overline{B} & c \\ \end{array} \right) | \, a,c \in {\mathbb Z} \, , \, B \in {\cal O}_{d} \, , \, ac - |B|^{2} < 0 \mbox{, and } (a,B,c) \mbox{ primitive } \right\}$$ 
Let $\Phi : {\cal H}_{d} \rightarrow \Sigma_{d}$ be the obvious $2$ to $1$ map.
We see that $det (A) = -{\cal D}(\Phi(A))$.
There is a natural action of $SL_{2}({\cal O}_{d})$ on ${\cal H}_{d}$ by $T \cdot A = T A T^{*}$.
As $\pm I$ are in the kernel of the action, we can induce an action by $PSL_{2}({\cal O}_{d})$.
This action clearly preserves determinants.

If $A \in {\cal H}_{d}$, a calculation shows that $\Phi(V \cdot A) = T \cdot \Phi(A)$, where $V = (T^{-1})^{t}$.
Thus, $\Phi$ descends to a map from the orbit space of ${\cal H}_{d}$ to the orbit space of $\Sigma_{d}$.
This in turn implies that the discriminant function on $\Sigma_{d}$ is invariant under the action of $PSL_{2}({\cal O}_{d})$.

Now we note that for any $T \in PSL_{2}({\cal O}_{d})$, $Stab_{PSL_{2}({\cal O}_{d})}(T \cdot {\cal C}) = T (Stab_{PSL_{2}({\cal O}_{d})}({\cal C})) T^{-1}$.
Since $Stab_{PSL_{2}({\cal O}_{d})}({\cal C})$ and $Stab_{PSL_{2}({\cal O}_{d})}({\cal C}')$ are commensurable if and only if they are conjugate, we obtain the following \cite{R1}.\\

\noindent
{\bf Theorem 4.2.1 } {\em Suppose $Stab_{PSL_{2}({\cal O}_{d})}({\cal C})$ and $Stab_{PSL_{2}({\cal O}_{d})}({\cal C'})$ are non-elementary Fuchsian groups commensurable in $PSL_{2}({\cal O}_{d})$.
Then ${\cal D}({\cal C}) = {\cal D}({\cal C'})$.}\\

Suppose $D \in {\mathbb Z}^{+}$ and that ${\cal C}_{D} \subset \widehat{\mathbb C}$ is a circle centered at the origin of ${\mathbb C}$ with radius $\sqrt{D}$.
Note that ${\cal D}({\cal C}_{D}) = D$.
One can check \cite{M} that
$$Stab_{PSL_{2}({\cal O}_{d})}({\cal C}_{D}) = \left\{ P \left( \begin{array}{cc}
\alpha & D \beta \\
\overline{\beta} & \overline{\alpha} \\ \end{array} \right)
| \, \alpha , \beta \in {\cal O}_{d} \mbox{ and } |\alpha|^{2} - D |\beta|^{2} = 1 \right\}$$
The following theorem provides the arithmetic structure of the groups $Stab_{PSL_{2}({\cal O}_{d})}({\cal C})$  (see \cite{M}, \cite{MR2}).\\

\noindent
{\bf Theorem 4.2.2 } {\em Suppose $d$ is a square-free integer.
Then, for every $D \in {\mathbb Z}^{+}$, $Stab_{PSL_{2}({\cal O}_{d})}({\cal C}_{D})$ is a Fuchsian group derived from the quaternion algebra
$$A = \left( \frac{-d,D}{\mathbb Q} \right)$$
Moreover, $Stab_{PSL_{2}({\cal O}_{d})}({\cal C}_{D})$ and $Stab_{PSL_{2}({\cal O}_{d})}({\cal C}_{D'})$ are commensurable in $PSL_{2}({\cal O}_{d})$ if and only if $D = D'$.}\\

\noindent
{\bf Proof: } Let ${\cal O} \subset A$ be the order given by 
$${\cal O} = \{ x = x_{0} + x_{1}{\bf i} + x_{2}{\bf j} + x_{3}{\bf k} \, | \, x_{0},x_{1},x_{2},x_{3} \in {\mathbb Z} \} \mbox{ if } d \equiv 1 \mbox{ or } 2 \mbox{ (mod 4) }$$
or
$${\cal O} = \{ x = \frac{x_{0}}{2} + \frac{x_{1}}{2}{\bf i} + \frac{x_{2}}{2}{\bf j} + \frac{x_{3}}{2}{\bf k} \, | \, x_{0},x_{1},x_{2},x_{3} \in {\mathbb Z} \mbox{ and } x_{0} \equiv x_{1} \, , \, x_{2} \equiv x_{3} \mbox{ (mod 2) } \} \mbox{ if } d \equiv 3 \mbox{ (mod 4) }$$

When $d \equiv 1,2$ (mod 4), ${\cal O}$ is easily seen to be an order in $A$.
When $d \equiv 3$ (mod 4), ${\cal O}$ is a finitely generated ${\mathbb Z}$-module and has ${\cal O} \otimes_{\mathbb Z} {\mathbb Q} = A$.  Further, one can check that $\rho({\cal O}) = R$, where
$$R = \left\{ \left( \begin{array}{cc}
\alpha & D \beta \\
\overline{\beta} & \overline{\alpha} \\ \end{array} \right)
| \, \alpha,\beta \in {\cal O}_{d} \right\}.$$
$R$ is a ring with $1$ since ${\cal O}_{d}$ is.
Therefore, ${\cal O}$ is a ring with $1$ and hence an order.

Now we see that $P \rho {\cal O}^{1} = Stab_{PSL_{2}({\cal O}_{d})}({\cal C}_{D})$ and the first assertions follows.
The second assertion is immediate from Theorem 4.2.1 and the fact that ${\cal D}({\cal C}_{D}) = D$. $\Box$\\

From Theorems 4.1.1 and 4.2.2 one can obtain the following useful criteria for $Stab_{PSL_{2}({\cal O}_{d})}({\cal C}_{D})$ to be co-compact (see \cite{MR1}).\\

\noindent
{\bf Theorem 4.2.3 } {\em Suppose $d,D \in {\mathbb Z}^{+}$ with $d \geq 3$ prime, and $D$ a quadratic non-residue (mod $d$).
Then $Stab_{PSL_{2}({\cal O}_{d})}({\cal C}_{D})$ is a co-compact Fuchsian group.}\\

\noindent
{\bf Remark: } Since the squaring endomorphism of $({\mathbb Z} / d {\mathbb Z})^{*}$ is not surjective when $d \geq 3$, quadratic non-residues always exist.

\subsection{Figure eight knot group}

It is well known \cite{Ri} that the complement of the figure eight knot in $S^{3}$ admits a complete, finite volume hyperbolic structure.
That is, $M_{8} \cong {\mathbb H}^{3} / \Gamma_{8}$, where $\Gamma_{8} \cong \pi_{1}(M_{8})$.
We make no distinction between $M_{8}$ and ${\mathbb H}^{3} / \Gamma_{8}$.
By conjugating if necessary, we may assume that $\Gamma_{8}$ is an index $12$ subgroup of $PSL_{2}({\cal O}_{3})$ and that
$$\Gamma_{8} = \left< P \left( \begin{array}{cc}
1 & 1 \\
0 & 1 \\ \end{array} \right),
P \left( \begin{array}{cc}
1 & 0 \\
-\omega & 1 \\ \end{array} \right) \right>$$
where $\omega^{2} + \omega + 1 = 0$ (note that ${\mathbb Z}[\omega] = {\cal O}_{3}$) \cite{Ri}.\\

\noindent
{\bf Theorem 4.3.1 } {\em For each positive integer $D \equiv 2$ (mod 3), $Stab_{\Gamma_{8}}({\cal C}_{D})$ is a co-compact Fuchsian group.
Moreover, $Stab_{\Gamma_{8}}({\cal C}_{D})$ and $Stab_{\Gamma_{8}}({\cal C}_{D'})$ are commensurable in $\Gamma_{8}$ if and only if $D = D'$.}\\

\noindent
{\bf Proof: } Let $D$ be as in the statement of the theorem.
Clearly $Stab_{\Gamma_{8}}({\cal C}_{D}) = Stab_{PSL_{2}({\cal O}_{3})}({\cal C}_{D}) \cap \Gamma_{8}$, so that
$$|Stab_{PSL_{2}({\cal O}_{3})}({\cal C}_{D}):Stab_{\Gamma_{8}}({\cal C}_{D})| \leq |PSL_{2}({\cal O}_{3}):\Gamma_{8}| = 12$$
Therefore, $Stab_{\Gamma_{8}}({\cal C}_{D})$ is co-compact if and only if $Stab_{PSL_{2}({\cal O}_{3})}({\cal C}_{D})$ is.
Since $2$ is not a square (mod 3), Theorem 4.2.3 implies $Stab_{PSL_{2}({\cal O}_{3})}({\cal C}_{D})$ is co-compact.
The second assertion follows from Theorem 4.2.2. $\Box$\\

We will need another fact concerning $\Gamma_{8}$ which is of an arithmetic nature.
Given any integer $n \geq 2$, let $R_{n}$ be the ring ${\cal O}_{3}/(n)$, where $(n)$ is the principal ideal in ${\cal O}_{3}$ generated by $n$.
For each $n \geq 2$, we define a homomorphism
$$\Phi_{n} : PSL_{2}({\cal O}_{3}) \rightarrow PSL_{2}(R_{n})$$
which is reduction of the entries modulo $(n)$.
The kernel of this homomorphism is a finite index, normal subgroup of $PSL_{2}({\cal O}_{3})$ called the {\em principal congruence subgroup of level n}, denoted $\Gamma(n)$.
Any subgroup of $PSL_{2}({\cal O}_{3})$ containing $\Gamma(n)$ for some $n \geq 2$ is called a {\em congruence subgroup}.
The theorem concerning the figure eight knot group which we need is the following well known fact.
We include a proof (this one due to Mark Baker) for completeness.\\

\noindent
{\bf Theorem 4.3.2 } {\em $\Gamma_{8}$ is a congruence subgroup of $PSL_{2}({\cal O}_{3})$ containing $\Gamma(4)$.}\\

\noindent
{\bf Proof: } Let $\Gamma_{8}' \subset PSL_{2}({\cal O}_{3})$ be the group
$$\Gamma_{8}' = \left< P \left( \begin{array}{cc}
1 & 1 \\
0 & 1 \\ \end{array} \right),
P \left( \begin{array}{cc}
1 & 0 \\
-\omega & 1\\ \end{array} \right),
P \left( \begin{array}{cc}
1 & 1+2 \omega \\
0 & 1\\ \end{array} \right) \right>$$

\noindent
According to \cite{Ri}, $\Gamma_{8} \subset \Gamma_{8}'$ is a subgroup of index $2$.
We prove the theorem by proving\\

{\bf{\em i}}. For any subgroup $H \subset \Gamma(2)$ with $|\Gamma(2):H| \leq 2$, we have $\Gamma(4) \subset H$,\\

{\bf{\em ii}}. $\Gamma(2) \subset \Gamma_{8}'.$\\\\
This will suffice since by ({\em ii}), $|\Gamma(2):\Gamma(2) \cap \Gamma_{8}| = |\Gamma(2) \cap \Gamma_{8}': \Gamma(2) \cap \Gamma_{8}| \leq |\Gamma_{8}':\Gamma_{8}| = 2$.
So by ({\em i}), $\Gamma(2) \cap \Gamma_{8}$ contains $\Gamma(4)$, so $\Gamma_{8}$ does.\\

\noindent
{\em Proof of (i): }  We prove this by showing that $\Gamma(4) = \Gamma(2)^{(2)}$, where $\Gamma(2)^{(2)} = < \gamma^{2} \, | \, \gamma \in \Gamma(2)>$.
For if this holds, then for any subgroup $H \subset \Gamma(2)$ of index no more than 2, we have $\gamma^{2} \in H$ for any $\gamma \in \Gamma(2)$.
Therefore, $\Gamma(4) = \Gamma(2)^{(2)} \subset H$ as required.

We first note that $G = \Gamma(2) / \Gamma(2)^{(2)} \cong ({\mathbb Z} / 2{\mathbb Z})^{n}$, for some $n \in {\mathbb Z}^{+} \cup \{ 0 \}$, since $\Gamma(2)$ is finitely generated.
A simple calculation shows that $\gamma^{2} \in \Gamma(4)$, $\forall \gamma \in \Gamma(2)$, so that $\Gamma(2)^{(2)} \subset \Gamma(4)$.
Hence, $\Gamma(2) / \Gamma(4)$ is a quotient of $G$.
From \cite{N}, we see that $|\Gamma(2) : \Gamma(4)| = 32$, so that $\Gamma(2) / \Gamma(4) \cong ({\mathbb Z} / 2{\mathbb Z})^{5}$, and $n \geq 5$.

In \cite{Bk1}, it is shown that $\Gamma(2) = \pi_{1}(S^{3} \setminus L)$ where $L$ is a five component link in $S^{3}$.
Therefore, $(\Gamma(2))^{ab} = {\mathbb Z}^{5}$.
Since $G$ is abelian, the quotient map $\Gamma(2) \rightarrow G$ factors through the abelianization.
Therefore, $n \leq 5$ implying $n=5$ and $\Gamma(2)^{(2)} = \Gamma(4)$.
Thus, ({\em i}) follows.\\

\noindent
{\em Proof of (ii): }  We consider 
$$K=core_{PSL_{2}({\cal O}_{3})}(\Gamma_{8}') = \bigcap_{\gamma \in PSL_{2}({\cal O}_{3})} \gamma \Gamma_{8}' \gamma^{-1} \lhd PSL_{2}({\cal O}_{3})$$
As $K \subset \Gamma_{8}'$, it will suffice to show that $\Gamma(2) \subset K$.
Note that $K = \bigcap_{i=1}^{6} s_{i} \Gamma_{8}' s_{i}^{-1}$ where $s_{1},...,s_{6}$ are coset representatives for $\Gamma_{8}'$ in $PSL_{2}({\cal O}_{3})$.
Moreover, both $PSL_{2}({\cal O}_{3})$ and $\Gamma_{8}'$ have a single conjugacy class of maximal elementary subgroups fixing a single point (the conjugacy class of the stabilizer of infinity).
Therefore, we can take $s_{1},...,s_{6}$ to lie in $Stab_{PSL_{2}({\cal O}_{3})}(\infty)$, that is, we may assume the $s_{i}$ are upper triangular.
Now, we also have
$$Stab_{\Gamma(2)}(\infty) = \left< P \left( \begin{array}{cc}
1 & 2 \\
0 & 1 \\ \end{array} \right),
P \left( \begin{array}{cc}
1 & 2 \omega \\
0 & 1 \\ \end{array} \right) \right>$$
It is immediate that $Stab_{\Gamma(2)}(\infty) \subset \Gamma_{8}'$.
Furthermore, normality of $\Gamma(2)$ in $PSL_{2}({\cal O}_{3})$ implies $s_{i} Stab_{\Gamma(2)}(\infty) s_{i}^{-1} = Stab_{\Gamma(2)}(\infty)$, hence $Stab_{\Gamma(2)}(\infty) \subset K$.

Now we again use that fact that $\Gamma(2) \cong \pi_{1}(S^{3} \setminus L)$.
A Wirtinger presentation generates $\pi_{1}(S^{3} \setminus L)$ by meridians \cite{Ro}, hence $\Gamma(2)$ is generated by parabolics.
If follows that the generators of $\Gamma(2)$ are $PSL_{2}({\cal O}_{3})$ conjugates of elements in $Stab_{\Gamma(2)}(\infty)$.
Since $K \lhd PSL_{2}({\cal O}_{3})$, these generators must lie in $K$.
Thus, $\Gamma(2) \subset K$, and ({\em ii}) follows, completing the proof.$\Box$

\subsection{Proof of Theorem 4.4.1}

Assuming Theorem 5.1.1, we prove\\

\noindent
{\bf Theorem 4.4.1 } {\em There exists a compact, orientable, irreducible $3$-manifold $M$, with torus boundary, having the following property.
Given any positive integer, $n$, there exists a closed, orientable, immersed, incompressible surface $F \looparrowright M$ with no incompressible annulus joining $F$ and $\partial M$, such that $F$ compresses in $M(\alpha)$ and $M(\beta)$ and $\Delta(\alpha,\beta) > n$.}\\
 
{\bf Proof: }  Let $M$ be the exterior of the figure eight knot in $S^{3}$, which is a compact, orientable, irreducible $3$-manifold with torus boundary.
$M_{8}$ is homeomorphic to the interior of $M$ and for the remainder of the proof, we regard $M_{8}$ as this subset of $M$.
Note that $M_{8}(\alpha) = M(\alpha)$ for any slope on $\partial M$.
Let $\mu$ and $\lambda$ denote the standard basis for $\pi_{1}(\partial M)$, so that slopes on $\partial M$ are represented by $\frac{p}{q} \in {\mathbb Q} \cup \{ \infty \}$. 

Now, let $p = 4$ and  $q > n$ be any odd integer.
By Theorem 5.1.1, there exists $F$, a closed, orientable, immersed totally geodesic surface in $M_{8}$ which compresses in $M_{8}(\frac{p}{q}) = M(\frac{p}{q})$.
Since $M_{8}(\frac{1}{0}) = S^{3}$ is simply connected, we have that $F$ compresses in $M_{8}(\frac{1}{0}) = M(\frac{1}{0})$.
We also have $\Delta(\frac{p}{q},\frac{1}{0}) = q > n$, so that to complete the proof, we need only prove that there exists no incompressible annulus in $M$ with one boundary component in $F$ and the other in $\partial M$.
Such an annulus defines a free homotopy from an essential curve in $F$ to an essential curve in $\partial M$.
That is, an essential curve in $F$ is peripheral.
This is forbidden by Corollary 3.3.2, hence no such annulus exists, and the theorem follows. $\Box$

\section{Compressions and normal closures}

\subsection{Compressions in $M_{8}(\alpha)$}

In this section we prove the main theorem.\\

\noindent
{\bf Theorem 5.1.1 } {\em Suppose $4|p$ and $3 \nmid p$.
Then, for any $q$ with $gcd(p,q)=1$ there exists infinitely many non-commensurable, closed, orientable, immersed, totally geodesic surfaces in $M_{8}$ which compress in $M_{8}(\frac{p}{q})$.}\\

\noindent
{\bf Remark: } We say that two surfaces in $M_{8}$ are commensurable if their fundamental groups are commensurable in $\Gamma_{8}$.\\

\noindent
{\bf Proof: } Let $\mu$ and $\lambda$ be elements in $Stab_{\Gamma_{8}}(\infty)$ representing the standard meridian-longitude basis for $\pi_{1}(\partial \overline{M_{8}})$.
It can be shown \cite{T}, that
$$ \mu = P \left( \begin{array}{cc}
1 & 1 \\
0 & 1 \\ \end{array} \right)
\hspace{1truein}
\lambda = P \left( \begin{array}{cc}
1 & 4 \omega + 2 \\
0 & 1 \\ \end{array} \right)$$

Let $p$ and $q$ be as in the statement of the theorem, and put
$$ \sigma = \mu^{p} \lambda^{q} = P \left( \begin{array}{cc}
1 & p + q (4 \omega + 2)\\
0 & 1 \\ \end{array} \right)
=
P \left( \begin{array}{cc}
1 & \xi\\
0 & 1 \\ \end{array} \right)$$
VanKampen's Theorem implies 
$$\pi_{1}(M_{8}(\frac{p}{q})) \cong \Gamma_{8} / \ll \sigma \gg$$
where $\ll \sigma \gg$ is the normal closure of $\{ \sigma \}$ in $\Gamma_{8}$.

The strategy is now the following.
We construct a sequence of distinct positive integers $D_{k} \equiv 2$ (mod $3$), and elements $g_{k} \in \ll \sigma \gg$ of the form
$$ g_{k} = P \left( \begin{array}{cc}
\alpha_{k} & D_{k} \beta_{k} \\
\overline{\beta_{k}} & \overline{\alpha_{k}} \\ \end{array} \right) \neq Id $$
Writing $ \Gamma_{D_{k}} = Stab_{\Gamma_{8}}({\cal C}_{D_{k}})$, Theorem 4.3.1 implies $\{ \Gamma_{D_{k}} \}$ is a sequence of non-commensurable (in $\Gamma_{8}$) co-compact Fuchsian subgroups of $\Gamma_{8}$.
Therefore, as in Sections 3.2 and 3.3, we obtain a sequence $\{ f_{k}: S_{\Gamma_{D_{k}}} \looparrowright M_{8} \}$ of pairwise non-commensurable, closed, orientable, immersed, totally geodesic surfaces in $M_{8}$ with $\psi \circ f_{k \, *} ( \pi_{1}(S_{\Gamma_{D_{k}}})) = \Gamma_{D_{k}}$.

The description of $Stab_{PSL_{2}({\cal O}_{3})}({\cal C}_{D})$ given in Section 4.2 and the fact that
$$\Gamma_{D_{k}} = \Gamma_{8} \cap Stab_{PSL_{2}({\cal O}_{3})}({\cal C}_{D_{k}})$$
together imply that $g_{k} \in \Gamma_{D_{k}}$.
So we see that $g_{k}$ represents a non-trivial element of $\pi_{1}(S_{\Gamma_{D_{k}}})$ lying in $\ll \sigma \gg$, and hence, $S_{D_{k}}$ compresses in $M_{8}(\frac{p}{q})$.\\

We now begin the construction of $D_{k}$ and $g_{k}$.
We have
$$|\xi|^{2} = (p + q (4 \omega +2))(p + q (4 \overline{\omega} +2)) = p^{2} + 12 q^{2}$$
Since $3 \nmid p$, we see that $3 \nmid |\xi|^{2}$, and $3 \nmid 4|\xi|^{2}$.
Therefore, $gcd(3,4|\xi|^{2})=1$ and there exists $r,t \in {\mathbb Z}$ such that
$$ -3r - 4|\xi|^{2}t=1$$
This implies that
$$ h = P \left( \begin{array}{cc}
\sqrt{-3} & 4 \xi t \\
\overline{\xi} & \sqrt{-3}r \\ \end{array} \right) $$
is an element of $PSL_{2}({\cal O}_{3})$.\\

\noindent
{\bf Claim: } $h \in \Gamma_{8}$.\\
To prove the claim, we consider the homomorphism $\Phi_{4}$ defining $\Gamma(4)$ (see Section 4.3).

By Theorem 4.3.2, we have that $\Gamma(4) \subset \Gamma_{8}$.
In particular $\Gamma_{8} = \Phi_{4}^{-1}(\Phi_{4}(\Gamma_{8}))$.
So to prove that $h \in \Gamma_{8}$, it suffices to prove that $\Phi_{4}(h) \in \Phi_{4}(\Gamma_{8})$.
We have
$$ h = P \left( \begin{array}{cc}
\sqrt{-3} & 4 \xi t\\
\overline{\xi} & \sqrt{-3}r\\ \end{array} \right)
\equiv
P \left( \begin{array}{cc}
\sqrt{-3} & 0 \\
2 & \sqrt{-3}\\ \end{array} \right)
\mbox{ (mod (4)) }$$
The congruence in the $(1,2)$-entry is clear.  The congruence in the $(2,1)$-entry follows from the fact that $4|p$ and $gcd(p,q)=1$ implies $q \equiv 1,3$ (mod 4), so that
$$ \overline{\xi} = p + q(4 \overline{\omega} + 2) \equiv 2q \equiv 2 \mbox{ (mod (4)) }$$
The congruence in the $(2,2)$-entry comes from
$$1 = det (h) = -3r - 4|\xi|^{2}t \equiv r \mbox{ (mod (4)) }$$

Next, we consider the following element $g \in \Gamma_{8}$, and its reduction modulo $(4)$
$$g = P \left( \left( \begin{array}{cc}
1 & 2\\
0 & 1\\ \end{array} \right)
\left( \begin{array}{cc}
1 & 0\\
\omega & 1\\ \end{array} \right) \right)^{2}
= P \left( \begin{array}{cc}
1+6\omega +4\omega^{2} & 4(1+\omega)\\
2(-1) & 1+2\omega\\ \end{array} \right)$$
$$= P \left( \begin{array}{cc}
\sqrt{-3}-4 & 2+2 \sqrt{-3} \\
-2 & \sqrt{-3}\\ \end{array} \right)
\equiv P \left( \begin{array}{cc}
\sqrt{-3} & 0\\
2 & \sqrt{-3}\\ \end{array} \right)
\mbox{ (mod (4)) }$$
Therefore, $\Phi_{4}(h) = \Phi_{4}(g)$, hence $\Phi_{4}(h) \in \Phi_{4}(\Gamma_{8})$ and the claim is established.\\

We now construct $D_{k}$ and $g_{k}$.
For each $k \in {\mathbb Z}^{+}$ set
$$ n_{k} = - 3 |\xi|^{2}(2+3k)+9$$
$$ D_{k} = |\xi|^{2} n_{k}^{2} + 2 +3k$$
$$g_{k} = \sigma^{n_{k}} (h \sigma h^{-1})^{6} \sigma^{n_{k}}$$
By construction, $g_{k} \in \ll \sigma \gg$ and $D_{k} \equiv 2$ (mod 3).
To complete the proof of the theorem, we must show that $g_{k}$ has the required form.

First, we compute
$$(h \sigma h^{-1})^{6} = h \sigma^{6} h^{-1} =
P \left( \left( \begin{array}{cc}
\sqrt{-3} & 4 \xi t\\
\overline{\xi} & \sqrt{-3} r\\ \end{array} \right)
\left( \begin{array}{cc}
1 & 6 \xi\\
0 & 1\\ \end{array} \right)
\left( \begin{array}{cc}
\sqrt{-3} r & -4 \xi t\\
-\overline{\xi} & \sqrt{-3}\\ \end{array} \right) \right)$$
$$= P \left( \left( \begin{array}{cc}
\sqrt{-3} & 6 \xi \sqrt{-3} + 4 \xi t\\
\overline{\xi} & 6 |\xi|^{2} + \sqrt{-3}r\\ \end{array} \right)
\left( \begin{array}{cc}
\sqrt{-3} r & -4 \xi t\\
-\overline{\xi} & \sqrt{-3}\\ \end{array} \right) \right)$$
$$= P \left( \begin{array}{cc}
-3r-\overline{\xi}(6 \xi \sqrt{-3} + 4 \xi t) & -4 \xi t \sqrt{-3} + \sqrt{-3}(6 \xi \sqrt{-3} + 4 \xi t)\\
 \overline{\xi} \sqrt{-3} r - \overline{\xi}(6 |\xi|^{2} + \sqrt{-3} r) & -4 |\xi|^{2} t + \sqrt{-3}(6 |\xi|^{2} + \sqrt{-3} r)\\ \end{array} \right)$$
$$= P \left( \begin{array}{cc}
-3r - 4 |\xi|^{2} t - 6 |\xi|^{2} \sqrt{-3} & -18 \xi \\
-6|\xi|^{2} \overline{\xi} & -4 |\xi|^{2}t -3r + 6 |\xi|^{2} \sqrt{-3} \\ \end{array} \right)$$
$$= P \left( \begin{array}{cc}
1 - 6|\xi|^{2} \sqrt{-3} & -18 \xi \\
-6|\xi|^{2}\overline{\xi} & 1+6|\xi|^{2}\sqrt{-3}\\ \end{array} \right)$$
This gives
$$g_{k} = \sigma^{n_{k}}( h \sigma h^{-1})^{6} \sigma^{n_{k}}$$
$$=
P \left( \left( \begin{array}{cc}
1 & n_{k} \xi\\
0 & 1\\ \end{array} \right)
\left( \begin{array}{cc}
1 - 6|\xi|^{2} \sqrt{-3} & -18 \xi \\
-6|\xi|^{2}\overline{\xi} & 1+6|\xi|^{2}\sqrt{-3}\\ \end{array} \right)
\left( \begin{array}{cc}
1 & n_{k} \xi\\
0 & 1 \\ \end{array} \right) \right)$$

$$=
P \left( \left( \begin{array}{cc}
1-6|\xi|^{2}\sqrt{-3} - 6 n_{k}|\xi|^{4} & -18 \xi + n_{k} \xi(1+6|\xi|^{2}\sqrt{-3})\\
-6|\xi|^{2}\overline{\xi} & 1 + 6 |\xi|^{2}\sqrt{-3}\\ \end{array} \right)
\left( \begin{array}{cc}
1 & n_{k}\xi\\
0 & 1\\ \end{array} \right) \right) $$

$$=
P \left( \begin{array}{cc}
1-6|\xi|^{2}\sqrt{-3} - 6 n_{k}|\xi|^{4} & n_{k}\xi( 1-6|\xi|^{2}\sqrt{-3} - 6 n_{k}|\xi|^{4} ) - 18 \xi + n_{k} \xi(1+6|\xi|^{2}\sqrt{-3})\\
-6|\xi|^{2}\overline{\xi} & - 6 n_{k}|\xi|^{4} + 1 + 6 |\xi|^{2}\sqrt{-3}\\ \end{array} \right)$$

$$= P \left( \begin{array}{cc}
1 - 6 n_{k}|\xi|^{4} - 6 |\xi|^{2}\sqrt{-3} & n_{k}\xi( 2 - 6 n_{k}|\xi|^{4} ) - 18 \xi\\
-6|\xi|^{2}\overline{\xi} & 1 - 6 n_{k}|\xi|^{4} + 6 |\xi|^{2}\sqrt{-3}\\ \end{array} \right)$$

Now set
$$ \alpha_{k} = 1 - 6 n_{k}|\xi|^{4} - 6 |\xi|^{2}\sqrt{-3}$$
$$ \beta_{k} = -6 |\xi|^{2}\xi$$
and note that
$$n_{k}\xi( 2 - 6 n_{k}|\xi|^{4} ) - 18 \xi = 2 n_{k} \xi - 6 n_{k}^{2}|\xi|^{4} \xi -18 \xi = (-3|\xi|^{2}(2+3k) + 9)2 \xi - |\xi|^{2} n_{k}^{2} \cdot 6 |\xi|^{2} \xi -18 \xi$$
$$ = (2 + 3k) (-6 |\xi|^{2} \xi) + 18 \xi + |\xi|^{2} n_{k}^{2}(- 6 |\xi|^{2} \xi) -18 \xi = (|\xi|^{2} n_{k}^{2} +2+3k)(-6 |\xi|^{2} \xi) = D_{k}\beta_{k}$$
Thus we have
$$ g_{k} = P \left( \begin{array}{cc}
1 - 6 n_{k}|\xi|^{4} - 6 |\xi|^{2}\sqrt{-3} & n_{k}\xi( 2 - 6 n_{k}|\xi|^{4} ) - 18 \xi\\
-6|\xi|^{2}\overline{\xi} & 1 - 6 n_{k}|\xi|^{4} + 6 |\xi|^{2}\sqrt{-3}\\ \end{array} \right)
= P \left( \begin{array}{cc}
\alpha_{k} & D_{k} \beta_{k}\\
\overline{\beta_{k}} & \overline{\alpha_{k}}\\ \end{array} \right)$$
This completes the proof.$\Box$\\

The element $h$ in the above proof was arrived at by attempting to ``match up'' an invariant hyperbolic plane of $\sigma$ with one from a conjugate of $\sigma$.
The group generated by $\sigma$ and $h \sigma h^{-1}$ is a non-elementary Fuchsian subgroup of $\Gamma_{8}$, which thus contains a rank-$2$ free subgroup of hyperbolic elements in $\ll \sigma \gg$ with real traces.
The fact that $\Gamma_{8}$ contains an abundance of co-compact Fuchsian subgroups makes it possible to construct nontrivial elements in the intersection of some of these groups with $\ll \sigma \gg$.

It seems worthwhile to compute $h$ and $g_{k}$ explicitly for some choice of $p,q$, and $k$.
An appendix is included at the end of the paper which contains $h$ and a list of $g_{k}$ and $D_{k}$ for $k \in \{ 1,...,10 \}$, when $p = 20$ and  $q = 7$.

\subsection{Related results}

The basic construction in the proof of Theorem 5.1.1 is quite general.
In particular, we have\\

\noindent
{\bf Theorem 5.2.1 } {\em Let $d \geq 3$ with $d$ prime.
Suppose also that
$$\sigma = P \left( \begin{array}{cc}
1 & \xi \\
0 & 1\\ \end{array} \right)
\in PSL_{2}({\cal O}_{d})$$
is such that $d \nmid |\xi|^{2}$.
Then $\ll \sigma \gg$ non-trivially intersects infinitely many non-commensurable (in $PSL_{2}({\cal O}_{d})$) co-compact Fuchsian subgroups of $PSL_{2}({\cal O}_{d})$.}\\

\noindent 
We only sketch the proof as it is almost identical to the proof of Theorem 5.1.1.\\

\noindent
{\bf Sketch of Proof: } As in the previous proof, we construct, for each $k \in {\mathbb Z}^{+}$, $g_{k} \in \ll \sigma \gg$ of the form
$$g_{k} = P \left( \begin{array}{cc}
\alpha_{k} & D_{k} \beta_{k} \\
\overline{\beta_{k}} & \overline{\alpha_{k}}\\ \end{array} \right)$$
where $\{ D_{k} \}$ is a sequence of distinct positive integers, each of which is a quadratic non-residue (mod d).
Therefore, $Stab_{PSL_{2}({\cal O}_{d})}({\cal C}_{D_{k}})$ is a co-compact Fuchsian group by Theorem 4.2.3, and $g_{k} \in \ll \sigma \gg \cap Stab_{PSL_{2}({\cal O}_{d})}({\cal C}_{D_{k}})$.
Moreover, the groups in the sequence $\{ Stab_{PSL_{2}({\cal O}_{d})}({\cal C}_{D_{k}}) \}$ are all non-commensurable in $PSL_{2}({\cal O}_{d})$ by Theorem 4.2.2, as required.

By hypothesis, there exists $r$ and $t$ such that $-d r - |\xi|^{2} t = 1$.
This implies that 
$$h = P \left( \begin{array}{cc}
\sqrt{-d} & \xi t \\
\overline{\xi} & \sqrt{-d} r \\ \end{array} \right) \in PSL_{2}({\cal O}_{d})$$

Let $x < d$ be a non-square (mod d) (see the remark at the end of Section 4.2).
For each $k \in {\mathbb Z}^{+}$ define
$$ n_{k} = - d |\xi|^{2}( d k + x ) + d^{2}$$
$$ D_{k} = n_{k}^{2}|\xi|^{2} + d k + x$$
$$ g_{k} = \sigma^{n_{k}} (h \sigma h^{-1})^{2d} \sigma^{n_{k}}$$
$$\alpha_{k} = 1 - 2 d n_{k}|\xi|^{4}- 2 d |\xi|^{2} \sqrt{-d}$$
$$\beta_{k} = - 2 d |\xi|^{2} \xi$$
By definition, $D_{k} \equiv x$ (mod d), and hence is not a square (mod d).
As in the previous proof, a calculation shows 
$$g_{k} = P \left( \begin{array}{cc}
\alpha_{k} & D_{k} \beta_{k} \\
\overline{\beta} & \overline{\alpha} \\ \end{array} \right)$$
thus completing the proof.$\Box$\\

Theorem 5.2.1 does not quite have the same topological implications as Theorem 5.1.1, since $M_{PSL_{2}({\cal O}_{d})}$ is never a manifold (i.e. $PSL_{2}({\cal O}_{d})$ always contains torsion).

We wish to use Theorem 5.2.1 to construct examples of cusped hyperbolic manifolds and arbitrarily large surgeries in which totally geodesic surfaces compress.
Selberg's Lemma (see \cite{Ra}, for example) guarantees the existence of a hyperbolic $3$-manifold $M_{\Gamma}$ that is a finite sheeted orbifold cover of $M_{PSL_{2}({\cal O}_{d})}$ (pass to $\Gamma \subset PSL_{2}({\cal O}_{d})$ a finite index, torsion free subgroup).
The problem in trying to use these covers is that the conjugating element $h$ in the proof of Theorem 5.2.1 need not lie in $\Gamma$.

There is a special situation of torsion free, finite index subgroups of $PSL_{2}({\cal O}_{d})$ where Theorem 5.2.1 can be applied.

Suppose $ d > 3 $ is prime.
In this situation, note that peripheral subgroup $Stab_{PSL_{2}({\cal O}_{d})}(\infty)$ is a rank-$2$ torsion free abelian group (the only units of ${\cal O}_{d}$ are $\pm 1$ for $d > 3$).
We say that a finite index subgroup $\Gamma \subset PSL_{2}({\cal O}_{d})$ is {\em $\infty$-non-separated} if we can choose coset representatives for $\Gamma \subset PSL_{2}({\cal O}_{d})$ from $Stab_{PSL_{2}({\cal O}_{d})}(\infty)$.
In this situation, we will denote such a set of coset representatives by $s_{1},...,s_{n}$.
This condition is equivalent to requiring that the index of $Stab_{\Gamma}(\infty)$ in $Stab_{PSL_{2}({\cal O}_{d})}(\infty)$ is equal to the index of $\Gamma$ in $PSL_{2}({\cal O}_{d})$.
This is also equivalent to the assumption that the preimage of the cusp of $M_{PSL_{2}({\cal O}_{d})}$ corresponding to $Stab_{PSL_{2}({\cal O}_{d})}(\infty)$ is a single cusp of $M_{\Gamma}$.

Suppose now that $ d > 3$ and $\Gamma \subset PSL_{2}({\cal O}_{d})$ is $\infty$-non-separated.
Note that any element $h \in PSL_{2}({\cal O}_{d})$ can be written as $g s_{i}$ for some $g \in \Gamma$ and some $i \in \{ 1,...,n \}$.
Let $\sigma$ be an arbitrary element of $Stab_{\Gamma}(\infty)$.
Since each $s_{i}$ centralizes $\sigma$, we see that
$$ h \sigma h^{-1} = g s_{i} \sigma s_{i}^{-1} g^{-1} = g \sigma g^{-1} $$
It follows that $\ll \sigma \gg_{\Gamma} = \ll \sigma \gg_{PSL_{2}({\cal O}_{d})}$ (here, $\ll \sigma \gg_{G}$ denotes the normal closure of $\sigma$ in $G$).\\

\noindent
{\bf Example: } It can be shown that there exist torsion free subgroups $\Gamma_{0}$ and $\Gamma_{1}$ in $PSL_{2}({\cal O}_{7})$, each of index $6$ (see \cite{GS} where these groups are called $\Gamma_{-7}(6,8)$ and $\Gamma_{-7}(6,9)$).
We show that for each $i=0,1$, there are infinitely many surgeries on $M_{\Gamma_{i}}$, such that for each such surgered manifold, $M_{\Gamma_{i}}(\alpha)$, there are infinitely many non-commensurable closed, orientable, immersed, totally geodesic surfaces in $M_{\Gamma_{i}}$ which compress in $M_{\Gamma_{i}}(\alpha)$.

Throughout, let $i=0,1$.
In \cite{GS} it is shown that $H_{1}(M_{\Gamma_{i}},{\mathbb Z}) \cong {\mathbb Z} \oplus ({\mathbb Z} / 3 {\mathbb Z})^{2}$.
A well known homological argument shows that a compact $3$-manifold $M$ with boundary must satisfy
$$b_{1}(M) \geq \frac{b_{1}(\partial M)}{2}.$$
It follows that $M_{\Gamma_{i}}$ can have no more than one cusp.
Since, $M_{\Gamma_{i}}$ must have at least one cusp (namely the one coming from $Stab_{\Gamma_{i}}(\infty)$), we see that $M_{\Gamma_{i}}$ has exactly one cusp, which implies $\Gamma_{i} \subset PSL_{2}({\cal O}_{7})$ is $\infty$-non-separated.

We take the basis for $Stab_{PSL_{2}({\cal O}_{d})}(\infty)$ given by
$$\mu = P \left( \begin{array}{cc}
1 & 1\\
0 & 1\\ \end{array} \right) \, , \,
\lambda = P \left( \begin{array}{cc}
1 & \eta\\
0 & 1\\ \end{array} \right) $$
where $\eta = \frac{ 1+\sqrt{-7}}{2}$.
Then any primitive element 
$$\sigma = \mu^{p} \lambda^{q} = P \left( \begin{array}{cc}
1 & \xi \\
0 & 1\\ \end{array} \right) \in Stab_{PSL_{2}({\cal O}_{d})}(\infty)$$
has $\xi = p + q \eta$ with $(p,q)$ a pair of coprime integers.
This implies
$$|\xi|^{2} = p^{2} + pq +2q^{2}$$
We assume that either $7 | p$ or $7 | q$ , so that $7 \nmid |\xi|^{2}$.

Taking $\sigma^{n_{i}}$ to be the smallest power of $\sigma$ which lifts to $Stab_{\Gamma_{i}}(\infty)$ we see that $n_{i}$ divides $6 = |Stab_{PSL_{2}({\cal O}_{d})}(\infty) : Stab_{\Gamma}(\infty)|$.
We can write
$$\sigma^{n_{i}} = \left( \begin{array}{cc}
1 & n_{i} \xi \\
0 & 1 \\ \end{array} \right) \mbox{ for } i= 0,1$$
Since $7 \nmid |\xi|^{2}$, we see that $7 \nmid |n_{i} \xi|^{2}$, so that $\sigma^{n_{i}}$ satisfies the hypothesis of Theorem 5.2.1 and hence $\ll \sigma^{n_{i}} \gg_{\Gamma_{i}} = \ll \sigma^{n_{i}} \gg_{PSL_{2}({\cal O}_{7})}$ intersects infinitely many non-commensurable co-compact Fuchsian subgroups.
As in the proof of Theorem 5.1.1, we see that this implies that there exist infinitely many non-commensurable closed, orientable, immersed, totally geodesic surfaces in $M_{\Gamma_{i}}$ which compress in $M_{\Gamma_{i}}(\sigma^{n_{i}})$.

This example also proves Theorem 4.4.1.
To see this, we refer to \cite{GS} for a presentation of $\Gamma_{0}$ ($=\Gamma_{-7}(6,8)$):
$$\Gamma_{0} = < g_{1},g_{2},g_{3} \, | \, g_{1}g_{2}^{-1}g_{3}g_{2}^{-1}g_{3}^{-1}g_{1}^{-1}g_{2}^{-1} \, , \, g_{1}g_{2}g_{3}^{-1}g_{1}g_{3}^{-1}g_{2}^{-1}g_{1}g_{3}^{-1} >$$
It is also shown there that $g_{2} = \mu^{2}$.
From this presentation, we see that $g_{2}$ is primitive in the abelianization, hence it is primitive in $\Gamma_{0}$.
Furthermore, a calculation shows $\Gamma_{0} / \ll g_{2} \gg \cong {\mathbb Z} * {\mathbb Z} / 3{\mathbb Z}$.
As this group can never contain the fundamental group of a surface (of positive genus), we see that every one of the surfaces found above compresses in $M_{\Gamma_{0}}(g_{2})$.
Therefore, if in the above construction we let $p=1$ and $q = 7k$ for $k \in {\mathbb Z}^{+}$, the sequence of slopes $(\sigma_{k})^{n_{0,k}} = (\mu \lambda^{7k})^{n_{0,k}}$ have totally geodesic surfaces compressing in $M_{\Gamma_{0}}(\sigma_{k}^{n_{0,k}})$, and 
$$7k = \Delta(\sigma_{k}, \mu) \leq \Delta(\sigma_{k}^{n_{0,k}},g_{2}) \rightarrow \infty \mbox{ as } k \rightarrow \infty.$$

\section{Appendix}

Here is a list of $g_{k}$ and $D_{k}$ from the proof of Theorem 5.1.1, for $p = 20$ and $q=7$ where $k \in \{ 1,...,10 \}$.
Note that it is not hard to see that these elements lie in $\Gamma(4) \subset \Gamma_{8}$.

First, we have
$$h = \left( \begin{array}{cc} 
\sqrt{-3}  &  -80-56 \sqrt{-3}\\
20-14\sqrt{-3} & 1317 \sqrt{-3}\\ \end{array} \right)$$

$$
\begin{array}{llc}
k         & D_{k}                   & g_{k}\\
\\

1 & 216733332353  & \left( \begin{array}{cc}
                    86746012705 - 5928 \sqrt{-3}   &  -25695903883771680 - 17987132718640176 \sqrt{-3}\\
                    -118560 + 82992 \sqrt{-3}      &  86746012705 + 5928 \sqrt{-3}\\ \end{array} \right)\\
\\

2 & 555090222500  & \left( \begin{array}{cc}
                    138825247393 - 5928 \sqrt{-3}  &  -65811496779600000 - 46068047745720000 \sqrt{-3} \\
                    -118560 + 82992 \sqrt{-3}      &  138825247393 + 5928 \sqrt{-3}  \\ \end{array} \right) \\
\\

3 & 1049684816711  & \left( \begin{array}{cc}
                     190904482081 - 5928 \sqrt{-3} &  -124450631869256160 - 87115442308479312 \sqrt{-3}\\
                     -118560 + 82992 \sqrt{-3}     &  190904482081 + 5928 \sqrt{-3}  \\ \end{array} \right) \\
\\

4 & 1700517114986  & \left( \begin{array}{cc}
                     242983716769 - 5928 \sqrt{-3} &  -201613309152740160 - 141129316406918112 \sqrt{-3}\\
                     -118560 + 82992 \sqrt{-3}     &  242983716769 + 5928 \sqrt{-3} \\ \end{array} \right) \\ 
\\

5 & 2507587117325  & \left( \begin{array}{cc}
                     295062951457 - 5928 \sqrt{-3} &  -297299528630052000 - 208109670041036400 \sqrt{-3}\\
                     -118560 + 82992 \sqrt{-3}     &  295062951457 + 5928 \sqrt{-3} \\ \end{array} \right) \\
\\

6 & 3470894823728  & \left( \begin{array}{cc}
                     347142186145 - 5928 \sqrt{-3} &  -411509290301191680 - 288056503210834176  \sqrt{-3}\\
                    -118560 + 82992 \sqrt{-3}      &  347142186145 + 5928 \sqrt{-3} \\ \end{array} \right)\\
\\

7 & 4590440234195  & \left( \begin{array}{cc}
                     399221420833 - 5928 \sqrt{-3} &  -544242594166159200 - 380969815916311440 \sqrt{-3}\\
                    -118560 + 82992 \sqrt{-3}      &  399221420833 + 5928 \sqrt{-3} \\ \end{array} \right)\\
\\

8 & 5866223348726 & \left( \begin{array}{cc}
                    451300655521 - 5928 \sqrt{-3}  &  -695499440224954560 - 486849608157468192 \sqrt{-3}\\
                    -118560 + 82992 \sqrt{-3}      &  451300655521 + 5928 \sqrt{-3} \\ \end{array} \right)\\
\\

9 & 7298244167321  & \left( \begin{array}{cc}
                     503379890209 - 5928 \sqrt{-3} &  -865279828477577760 - 605695879934304432 \sqrt{-3}\\
                    -118560 + 82992 \sqrt{-3}      &  503379890209 + 5928 \sqrt{-3} \\ \end{array} \right)\\
\\

10 & 8886502689980 & \left( \begin{array}{cc}
                     555459124897 - 5928 \sqrt{-3} &  -1053583758924028800 - 737508631246820160 \sqrt{-3}\\
                    -118560 + 82992 \sqrt{-3}      &  555459124897 + 5928 \sqrt{-3} \\ \end{array} \right)\\ \end{array}$$
\\

\noindent
Address:\\
Department of Mathematics,\\
University of Texas\\
Austin, TX 78712\\
Phone:\\
(512)-475-9146\\
email:\\
clein@math.utexas.edu\\


\begin{thebibliography}{99999}




\bibitem{Bk1} M. Baker, {\em Link complements and imaginary quadratic number fields}, preprint.\\
\bibitem{Bt} A. Bart, {\em Surface groups in surgered manifolds}, to appear in Topology.\\
\bibitem{BP} R. Benedetti and C. Petronio, {\em Lectures on Hyperbolic Geometry}, Springer-Verlag, Berlin Heidelberg (1992).\\
\bibitem{CGLS} M. Culler, C. Gordon, J. Luecke, and P. B. Shalen, {\em Dehn surgery on knots}, Annals of Math. {\bf 125} (1987), 237--300. \\
\bibitem{GS} F. Grunewald and J. Schwermer, {\em Subgroups of Bianchi groups and arithmetic quotients of hyperbolic $3$-space}, Trans. Am. Math. Soc. {\bf 335} (1993), 47--78.\\
\bibitem{H} J. Hempel, {\em $3$-Manifolds}, Princeton University Press, Princeton, New Jersey (1976).\\
\bibitem{M} C. Maclachlan, {\em Fuchsian subgroups of the groups $PSL_{2}({\cal O}_{d})$, in Low Dimensional Topology and Kleinian Groups}, ed D.B.A. Epstein LMS Lecture Notes Series {\bf 112} (1986), 305--311.\\
\bibitem{MR1} C. Maclachlan and A. W. Reid, {\em The Arithmetic of Hyperbolic $3$-Manifolds}, to appear Springer-Verlag.\\
\bibitem{MR2} C. Maclachlan and A. W. Reid, {\em Parameterizing Fuchsian subgroups of the Bianchi groups}, Can. J. Math. {\bf 1} (1991), 158--181.\\
\bibitem{Ma} A. Marden, {\em The geometry of finitely generated Kleinian groups}, Annals of Math. {\bf 99} (1974), 383--461.\\
\bibitem{N} M. Newmann, {\em Integral Matrices}, Academic Press, New York and London (1972). \\ 
\bibitem{Ra} J. Ratcliffe, {\em Foundations of Hyperbolic Manifolds}, Springer-Verlag, New York (1994).\\
\bibitem{R1} A. W. Reid, Ph.D. thesis, Univ. Aberdeen (1987).\\
\bibitem{R2} A. W. Reid, {\em Totally geodesic surfaces in hyperbolic $3$-manifolds}, Proc. Edinburgh Math. Soc. {\bf 34} (1991), 77--88.\\
\bibitem{Ri} R. Riley, {\em A quadratic parabolic group}, Math. Proc. Camb. Phil. Soc. {\bf 77} (1975), 281-288.\\
\bibitem{Ro} D. Rolfsen, {\em Knots and Links}, Publish or Perish (1977).\\
\bibitem{T} W.P. Thurston, {\em The Geometry and Topology of 3-manifolds}, Princeton University mimeo\-graphed notes (1979).\\
\bibitem{Wu} Y. Q. Wu, {\em Incompressibility of surfaces in surgered $3$-manifolds}, Topology {\bf 31} (1992), 271--279.\\




\end{thebibliography}
\end{document}